\documentclass[fleqn]{mat01}
\usepackage{times,mathtimy,amssymb}
\begin{document}

\setcounter{page}{199}
\firstpage{199}

\font\zz=msam10 at 10pt
\def\Box{\mbox{\zz{\char'244}}}

\renewcommand{\theequation}{\thesection\arabic{equation}}

\newcommand{\lr}{\longrightarrow}
\newcommand{\ra}{\rightarrow}
\newcommand{\hra}{\hookrightarrow}

\newtheorem{lem}{Lemma}
\newtheorem{thore}{Theorem}
\newtheorem{propo}{\rm PROPOSITION}
\newtheorem{coro}{\rm COROLLARY}
\newtheorem{theor}{\bf Theorem}

\font\zzzz=tibi at 13.5pt
\def\c{\mbox{\zzzz{c}}}
\def\E{\mbox{\zzzz{E}}}
\def\F{\mbox{\zzzz{F}}}

\font\ccc=rmtmib at 13.5pt
\def\ella{\mbox{\ccc{\char'140}}}
\def\Pia{\mbox{\ccc{\char'65}}}

\title{Containment of $\c_{\bf 0}$ and $\ell_{\bf 1}$ in $\Pi_{\bf 1}
\hbox{\bf (}\E\hbox{\bf ,}\ \F\hbox{\bf )}$}

\markboth{Mohsen Alimohammady}{Containment of $c_{0}$ and $\ell_{1}$ in
$\Pi_{1} (E, F)$}

\author{MOHSEN ALIMOHAMMADY}

\address{Department of Mathematics,
University of Mazandran, Faculty of Basic Sciences, Babolsar, Iran\\
\noindent E-mail: amohsen@umz.ac.ir}

\volume{114}

\mon{May}

\parts{2}

\Date{MS received 21 April 2004}

\begin{abstract}
Suppose $\Pi_{1} (E, F)$ is the space of all absolutely 1-summing
operators between two Banach spaces $E$ and $F$. We show that if $F$ has
a copy of $c_{0}$, then $\Pi_{1} (E, F)$ will have a copy of $c_{0}$,
and under some conditions if $E$ has a copy of $\ell_{1}$ then $\Pi_{1}
(E, F)$ would have a complemented copy of $\ell_{1}$.
\end{abstract}

\keyword{Absolutely 1-suuming operators; copy of $c_{0}$; copy of
$\ell_{1}$.}

\maketitle

\section{Introduction}

Many studies on copy of $c_{0}$ or $\ell_{1}$ in spaces of bounded
operators, weakly compact operators, and compact operators have been
made \cite{2,3,4,5,7}. Here we intend to prove similar results in the
space of absolutely 1-summing operators.

Absolutely 1-summing operators were introduced and studied by
Grothendi\'{e}ck in his famous r\'{e}sume \cite{6}. For two Banach
spaces $E$ and $F$, the operator $T\hbox{:}\ E \ra F$ is called
absolutely 1-summing operator, if given any sequence $(x_{n})_{n}$ of
$E$ for which $\sum_{n = 1}^{\infty} |x^{*} x_{n}| < \infty$, for each
$x^{*} \in E^{*}$, we have $\sum_{n = 1}^{\infty} \left\Vert Tx_{n}
\right\Vert < \infty.\ \Pi_{1} (E, F)$ is the Banach space of all
absolutely 1-summing operators from $E$ to $F$, endowed with
$\pi_{1}$-norm as follows:
\begin{align*}
\pi_{1} (T) &= {\inf} \ \left\lbrace \rho > 0 : \sum\limits_{i = 1}^{n}
\left\Vert T x_{i}\right\Vert\right.\\[.2pc]
&\leq \left. \rho \ {\sup}_{x^{*} \in B^{*}_{E}} \
\left(\sum\limits_{i = 1}^{n}|x^{*}(x_{i})|\right); \ x_{i}\in E, 1 \leq
i \leq n\right\rbrace.
\end{align*}
By Grothendi\'{e}ck--Pietsch Dominated Theorem \cite{1}, for $T \in
\Pi_{1} (E, F)$ there exists a regular probability measure $\mu$
defined on $B_{E^{*}}$ (with its weak$^{*}$ topology) for which
\begin{equation*}
\left\Vert T(x)\right\Vert\leq \pi_{1} (T) \int_{B_{E^{*}}} |x^{*} x| \
{\rm d}\mu (x^{*}),
\end{equation*}
holds for each $x \in E$. Therefore, $\left\Vert T(x) \right\Vert\leq
\pi_{1} (T) \left\Vert  x \right\Vert$.

The following lemma discusses the main results of our study.

\begin{lem}\hskip -.3pc {\rm \cite{1}}.\ \ Let $x^{*}$ and $y$ be{\rm ,}
in $B_{E^{*}}$ and $B_{F}${\rm ,} the unit balls of two Banach spaces
$E^{*}$ and $F$ respectively. Then $x^{*} \otimes y\in \Pi_{1} (E, F)$
and $\pi_{1} (x^{*} \otimes y) \leq 1$.\hfill$\Box\ $
\end{lem}

\section{Main results}

Here, we present a result showing complemented copy of $\ell_{1}$ in
$\Pi_{1} (E, F)$.

\setcounter{theor}{1}
\begin{theor}[\!]
Suppose $(x^{*}_{n})_{n}$ is a sequence in $B_{E^{*}}$ such that
$x_{n}^{*} x_{n} = 1$ and $x^{**} x_{n}^{*} > \epsilon$ for each $n \in
N${\rm ,} where $(x_{n})_{n}$ is a weakly unconditionally Cauchy
sequence in $E$ and $x^{**} \in E^{**}$ {\rm (}for example $E = c_{0},
e_{n} = x_{n}, x_{n}^{*} = e_{n} \in \ell_{1}, x^{**} = \chi N \in
\ell_{\infty})$. Then $\Pi_{1} (E, F)$ has a complemented copy of
$\ell_{1}$ if $F$ has a copy of $\ell_{1}$.
\end{theor}

\begin{proof}
From the assumption there is an $\ell_{1}$-basic sequence $(y_{n})_{n}$
in $F$ and so $M_{1}, M_{2} > 0$ such that for each $(a_{n})_{n} \in
\ell_{1}$ we have
\begin{equation*}
M_{1} \sum\limits_{n = 1}^{\infty} |a_{n}| \leq \left\Vert\sum\limits_{n
= 1}^{\infty} a_{n} y_{n}\right\Vert \leq M_{2} \sum\limits_{n =
1}^{\infty} |a_{n}|.
\end{equation*}
Therefore,
\begin{align*}
\pi_{1} \left(\sum_{n = 1}^{m} a_{n} x_{n}^{*} \otimes y_{n}\right)
&\leq \sum\limits_{n = 1}^{m} |a_{n}| \pi_{1} (x^{*}_{n} \otimes
y_{n})\\[.2pc]
&\leq \sum\limits_{n = 1}^{m} |a_{n}| \Vert x_{n}^{*} \Vert \cdot
\left\Vert y_{n} \right\Vert \leq M_{2} \sum\limits_{n = 1}^{m} |a_{n}|.
\end{align*}
On the other hand,
\begin{align*}
\pi_{1} \left(\sum\limits_{n = 1}^{m} a_{n} x_{n}^{*} \otimes
y_{n}\right) &\geq \left\Vert\sum\limits_{n = 1}^{m} x_{n}^{*} \otimes
y_{n}\right\Vert\\[.2pc]
&= {\sup}_{x^{**} \in B_{E^{**}}} \left\Vert \sum\limits_{n = 1}^{m}
a_{n} x^{**} (x_{n}^{*}) y_{n}\right\Vert\\[.2pc]
&\geq M_{1} \sum\limits_{n = 1}^{m} |a_{n}| \cdot |x_{0}^{**}
x_{n}^{*}| \geq CM_{1} \sum\limits_{n = 1}^{m} |a_{n}|.
\end{align*}

This shows that $\Pi_{1} (E, F)$ contains $(x^{*}_{n} \otimes
y_{n})_{n}$ as a copy of $\ell_{1}$. We show this copy is complemented
in $\Pi_{1}(E, F)$. Define $P\hbox{:}\ \Pi_{1} (E, F) \rightarrow [x_{n}^{*}
\otimes y_{n}]$ by $P(T) = \sum_{n = 1}^{\infty} (x_{n}^{*} \otimes
y_{n}) y_{n}^{*} T(x_{n})$, where $[x_{n}^{*} \otimes y_{n}]$ is the
closed linear span of $(x^{*}_{n} \otimes y_{n})$ in $\Pi_{1} (E, F)$.
Since $(x_{n})_{n}$ is $w u C$ and $T$ is the absolutely 1-summing
operator, $P$ is a well-defined linear map. From the closed graph
theorem one can easily show that $P$ is a projection from $\Pi_{1} (E,
F)$ onto $[x_{n}^{*} \otimes y_{n}]$ which completes the proof. \hfill
$\Box$
\end{proof}

Now we explain this result for asymptotically isomorphic copy of
$\ell_{1}$.\vspace{.5pc}

\setcounter{propo}{2}
\begin{propo}$\left.\right.$\vspace{.5pc}

\noindent $\Pi (E, F)$ has an asymptotically isomorphic copy of
$\ell_{1}$ if $F$ has an asymptotically isomorphic copy of $\ell_{1}$
too.
\end{propo}

\begin{proof}
Suppose for $(y_{n})$ in $F$ there is a positive null sequence
$(\epsilon_{n})$ such that
\begin{equation*}
\sum (1 - \epsilon_{n}) |a_{n}| \leq \left\Vert\sum a_{n}
y_{n}\right\Vert\leq \sum |a_{n}|.
\end{equation*}
Consider $x_{0}^{*}$ and $x_{0}$ in the unit sphere of $E^{*}$ and $E$
respectively such that $x_{0}^{*} x_{0} = 1$. Then
\begin{equation*}
\sum (1 - \epsilon_{n}) |a_{n}| \leq \left\Vert \sum a_{n} x_{0}^{*} x_{0}
y_{n}\right\Vert \leq \pi_{1} \left(\sum a_{n} x_{0}^{*} \otimes y_{n} \right) \leq
\sum |a_{n}|
\end{equation*}
which completes the proof. \hfill $\Box$\vspace{.4pc}
\end{proof}

\setcounter{theor}{3}
\begin{theor}[\!]
Suppose $L (E, F)$ contains a sequence $(T_{n})$ equivalent to the
standard unit vector basis of $c_{0}$ such that for $x_{0} \in E$
{\rm (}respectively $y_{0}^{*} \in F^{*}) T_{n} x_{0}$ {\rm (}respectively
$T_{n}^{*} y_{0}^{*}${\rm )} are basic sequences. Then $(F)$ or $E^{**}$
has a copy of $c_{0}$.
\end{theor}

\begin{proof}
Suppose $(T_{n} x_{0})$ is a basic sequence and $(y_{n}^{*})$ is its
coefficient functional. We can assume $\Vert
y_{n}^{*} \Vert\leq M$. On the other hand,
\begin{equation*}
C_{1} \hbox{sup}_{n} |a_{n}| \leq \left\Vert \sum\limits_{n + 1}^{\infty}
a_{n} T_{n}\right\Vert \leq C_{2} \hbox{sup}_{n} |a_{n}|.
\end{equation*}
Therefore,
\begin{equation*}
(C_{1}/M) |a_{n}| = (1/M) \bigg\vert y_{n}^{*} \left(\sum\limits_{n =
1}^{\infty} a_{n} T_{n}(x_{0})\right)\bigg\vert\leq
\left\Vert\sum\limits_{n = 1}^{\infty} a_{n} T_{n}\right\Vert \leq C\! -\! 2
\hbox{sup}_{n} |a_{n}|.
\end{equation*}

$\left.\right.$\vspace{-1.5pc}

\noindent This shows that $(T_{n} x_{0})$ is also a copy of $c_{0}$.
\hfill $\Box$\vspace{.4pc}
\end{proof}

\setcounter{coro}{4}
\begin{coro}$\left.\right.$\vspace{.5pc}

\noindent Let $({T_{n}})$ be a copy of $c_{0}$ in $\Pi(E, F)$ and
$x_{0}\in B_{E}$ such that $T_{n}(x_{0})$ is a semi-normalized sequence
in $F$. Then $F$ would have a copy of $c_{0}$.
\end{coro}

\begin{proof}
$M_{1}, M_{2} > 0$ such that for any $(a_{n})\in c_{0}$,
\begin{equation*}
M_{1}\ \hbox{sup}_{n}|a_{n}|\leq\pi_{1}\left(\sum_{n =
1}^{\infty}a_{n}T_{n}\leq M_{2}\ \hbox{sup}_{n}|a_{n}|\right).
\end{equation*}
But it is easy to see that $T_{n}(x_{0})$ is weakly null. From the
assumption $\lim\inf \|T_{n}x_{0}\| = C > 0$, it follows from
Bessaga--Pe\'{l}czynski' Theorem \cite{1} that it would be a basic
sequence with coefficient functionals $(y_{n}^{*})$ such that
$\|y_{n}^{*}\| < 1/C$. We would have,
\begin{align*}
|a_{n}| = \bigg\vert y_{n}^{*} \left( \sum\limits_{m = 1}^{\infty} a_{m}
T_{m} (x_{0}) \right) \bigg\vert &\leq (1/C) \left\Vert \sum\limits_{n =
1}^{\infty} a_{m} T_{m} (x_{0}) \right\Vert\\[.3pc]
&\leq (1/C) \left\Vert \sum\limits_{n = 1}^{\infty} a_{m} T_{m}
\right\Vert\\[.3pc]
&\leq (1/C) \pi_{i} \left(\sum_{n = 1}^{\infty} a_{m} T_{m} \right)\\[.3pc]
&\leq M - (2/C)\ \hbox{sup}_{n} |a_{n}|,
\end{align*}
which shows that $(T_{n}(x_{0}))$ is a $c_{0}$-basic sequence in
$F$.\hfill$\Box$
\end{proof}

\setcounter{theor}{5}
\begin{theor}[\!]
$\Pi_{1}(E, F)$ has a complemented copy of $c_{0}$ if $F$ has a
complemented copy of~$c_{0}$.
\end{theor}

\begin{proof}
Similar to the proof of Theorem~2 we may assume $(y_{n})_{n}$ is a
complemented copy of $c_{0}$ in $F$ and $(y_{n}^{*})_{n}$ the
coefficient functional of it in $[y_{n}]$ the closed linear span of
$(y_{n}^{*})_{n}$. Since $[y_{n}]$ is complemented in $F$, we can assume
$(y_{n}^{*})$ is a weak$^{*}$ null convergence sequence in $F^{*}$.
Define $P\hbox{:}\ \Pi_{1}(E, F)\rightarrow [x^{*} \otimes y_{n}]$ by
$P(T) (x) = \sum_{n = 1}^{\infty}(x^{*} \otimes y_{n})y_{n}^{*}T(x)$,
where $x^{*}$ is an arbitrary element in the unit sphere of $E^{*}$ and
$[x^{*}\otimes y_{n}]$ is the closed linear subspace generated by
$(x^{*}\otimes y_{n})$ in $\Pi_{1}(E, F)$. For any $(x_{i})_{i = 1}^{n}$
in $E$, we have
\begin{align*}
\sum\limits_{i = 1}^{n} \|PT (x_{i})\| &= \sum\limits_{i = 1}^{n}
\bigg\| \sum\limits_{n = 1}^{\infty} (y_{n}^{*} T(x_{i})) x^{*} \otimes
y_{n} \bigg\|\\[.2pc]
&\leq \sum\limits_{i = 1}^{n} C_{2} \hbox{sup}_{n}| (y_{n}^{*} T(x_{i}))|\\[.2pc]
&= C_{2} \sum\limits_{i = 1}^{n} \| T(x_{i})\| \leq C_{2} \pi_{1}(T)
\sum\limits_{i = 1}^{n} |x^{*}x_{i}|.
\end{align*}

This shows that $\pi_{1}(PT) \leq C_{2}\pi_{1}(T)$. Therefore, $P$ is a
bounded projection, which completes the proof.\hfill$\Box$
\end{proof}

\end{document}